\documentclass[11pt,notitlepage,twoside,a4paper]{amsart}
 \usepackage{amsfonts}

\usepackage{amsmath,amssymb,enumerate}

\usepackage{epsfig,fancyhdr,color}

\usepackage{amssymb}
\usepackage{amsmath,amsthm}
\usepackage{latexsym}
\usepackage{amscd}
\usepackage{psfrag}
\usepackage{graphicx}
\usepackage[latin1]{inputenc}
\usepackage[all]{xy}
\usepackage[mathcal]{eucal}

\definecolor{NoteColor}{rgb}{1,0,0}

\renewcommand{\textsc}{\textcolor{red}}

\newtheorem{theorem}{\rm\bf Theorem}[section]

\newtheorem*{theorem 1}{\rm\bf Proposition 1}
\newtheorem*{theorem 2}{\rm\bf Proposition 2}

\theoremstyle{definition}

\theoremstyle{remark}

\def\interieur#1{\mathord{\mathop{\kern 0pt #1}\limits^\circ}}

\title[On Lobachevsky's trigonometric formulae]{On Lobachevsky's trigonometric formulae}

\author{Athanase Papadopoulos}
\address{Athanase Papadopoulos,  Institut de Recherche Math\'ematique Avanc\'ee,
Universit\'e de Strasbourg and CNRS,
7 rue Ren\'e Descartes,
 67084 Strasbourg Cedex, France.} 
 \email{athanase.papadopoulos@math.unistra.fr}

\date{\today}

\begin{document}

\begin{abstract}  We elaborate on some important ideas contained in Loba\-chevsky's \emph{Pangeometry} and in some of his other memoirs. The ideas include the following: (1) The trigonometric formulae, which express the dependence between angles and edges of triangles, are not only tools, but they are used as the basic elements of any geometry. In fact, Lobachevsky developed a large set of analytical and geometrical theorems in non-Euclidean geometry using these formulae. (2) Differential and integral calculus are developed in hyperbolic space without the use of any Euclidean model of hyperbolic space. (3)  There exist models of spherical and of Euclidean geometry within hyperbolic geometry, and these models are used to prove the hyperbolic trigonometry formulae. (4)  If hyperbolic geometry were contradictory, then either Euclidean or spherical geometry would be contradictory. 
\\
We shall also see that some of these ideas were rediscovered by later mathematicians. 
\noindent  

\bigskip 

\noindent AMS classification: 01-99 ; 53-02 ; 53-03 ; 53A35.

\bigskip 

\noindent Keywords: Lobachevsky geometry ; non-Euclidean geometry ; hyperbolic geometry ; non-Euclidean trigonometric formulae ; models of hyperbolic space ; horosphere ; Lobachevsky integrals ; non-contradiction of geometry.

\bigskip 

\noindent  \emph{Acknowledgements.---} The author is partially supported by the French ANR project FINSLER. He wishes to thank the organizers of the International Seminar on History of Mathematics held in Delhi on November 19-20, 2012 and he is is most grateful to Norbert A'Campo who taught him non-Euclidean geometry.

\bigskip 

\noindent  The paper in final form will appear in \emph{Ga\d{n}ita Bh\=ar\=at\=\i }, the  Bulletin of the Indian Society for History of Mathematics.

 \end{abstract}
\maketitle

\section{Introduction}
Among the three founders of hyperbolic geometry - Lobachevsky, Bolyai and Gau\ss \  -- Lobachevsky was the first to publish a complete treatise on the subject, namely, his \emph{Elements of Geometry} \cite{Loba-Elements} (1829). In fact, three years before that publication (more precisely, on February 23\footnote{This date is according to the Gregorian calendar. The reader should be aware of the fact that concerning Lobachevsky, some confusion in dates may occur because in nineteenth century Russia, the Julian calendar was used, and in mentioning dates, some sources use the Julian calendar and others the Gregorian calendar.}, 1826), Lobachevsky gave a lecture (in French) at the Physical and  Mathematical Section of Kazan University, titled {\it Exposition succinte des principes de la g\'eom\'etrie avec une d\'emonstration rigoureuse du th\'eor\`eme des parall\`eles} (A Brief Exposition of the Principles of Geometry with a Rigorous Proof of the Theorem on Parallels). The exact content of the lecture is lost, but we know from Lobachevsky's own report in his \emph{Elements of Geometry} that the first part of that treatise is based on the material of the 1826 lecture.\footnote{On p. 1 of the {\it Elements of Geometry} \cite{Loba-Elements} Lobachevsky writes that this work is ``Extracted by the author himself from a paper  which he read on February 12, 1826, at a meeting of the Section for Physico-Mathematical Sciences, with the title \emph{Exposition succinte des principes de la G\'eom\'etrie, avec une d\'emonstration rigoureuse du th\'eor\`eme des parall\`eles}."
Likewise, in  the Introduction to the {\it New Elements of Geometry} \cite{Loba-New}  (1836), Lobachevski writes the following:
``Believing myself to have completely solved the difficult question, I wrote a paper on it in the year 1826: \emph{Exposition succinte des principes de la G\'eom\'etrie, avec une d\'emonstration rigoureuse du th\'eor\`eme des parall\`eles}, read on February 12, 1826, in the s\'eance of the Physico-mathematics faculty of the University of Kazan but nowhere printed."
} Lobachevsky's \emph{Elements of Geometry}, together with his later writings on hyperbolic geometry, constitute  a corpus of memoirs which is much more than an introduction to hyperbolic geometry. They contain a complete development of differential and integral calculus in a non-Euclidean setting and they reflect their author's profound vision on mathematics. We shall dwell on that in this article. Let us mention right away Gau\ss 's words, in his letter to the astronomer H. C. Schumacher, dated 28 November 1846, talking about one of Lobachevsky's memoirs\footnote{The memoir in question is Lobachevsky's \emph{Geometrische
Untersuchungen zur Theorie der Parallellinien} \cite{Loba-Geometrische}.}:  \begin{quote}\small
In developing the subject, the author followed a road different from the one I took myself; Lobachevsky carried out the task in a masterly fashion and in a truly geometric spirit. I consider it a duty to call your attention to this book, since I have no doubt that it will give you a tremendous pleasure ...\footnote{[Materiell f\"ur mich Neues habe ich also im Lobatschewskyschen Werke nicht gefunden, aber die Entwickelung ist auf anderm Wege gemacht, als ich selbst eingeschlagen habe, und zwar von Lobatschewsky auf eine meisterhafte Art in \"acht geometrischem Geiste. Ich glaube Sie auf das Buch aufmerksam machen zu m\"ussen, welches Ihnen gewiss ganz exquisiten Genuss gew\"ahren wird ...] Gaus's correspondence is published in Volume VIII of his \emph{Collected Works} \cite{Gauss}.}\end{quote}

During his lifetime, Lobachevsky never saw his work acknowledged and, in fact, nobody read his work seriously (with the notable exception of Gau\ss , who read some of it but remained silent about it, except for a private correspondence with his friends Schumacher and Taurinus; Lobachevsky was not aware of the fact that Gau\ss \ read his memoir). The reason was certainly the strong belief within the mathematical community that Euclid's parallel axiom is a consequence of the other axioms of Euclidean geometry and, therefore, that a geometry in which one takes as an axiom the negation of the parallel axiom, leaving all the other axioms of Euclidean geometry (that is, the Lobachevsky geometry), would be self-contradictory.\footnote{``Proofs" of the parallel axiom were given by Legendre, Bertrand and other important mathematicians, and there were several attempts, by Euler, Lagrange, Fourier and others; see \cite{Pont} and the commentary in \cite{L}.}

After his \emph{Elements of Geometry} \cite{Loba-Elements}, Lobachevsky wrote several memoirs on the same subject, reworking some of the proofs, improving some details,  and highlighting applications of his work outside the world of pure geometry,\footnote{For instance, Lobachevsky carried out a large amount of computations of areas and volumes of figures in hyperbolic space that lead to some new identities between integrals.  Integral identities were fashionable at that time. Despite this fact, Lobachevsky's work did not attract much attention. Some of his integration formulae that concern volumes of hyperbolic polyhedra were revived in the work of Coxeter, in the 1930s,  cf. \cite{Coxeter-Functions}, and they were used in the 1970s in the work of Thurston, see in particular Chapter 7 of \cite{Thurston}; see also Milnor \cite{Milnor}.} with the hope of attracting the attention of his colleagues. His efforts in that direction were vain, and his work was acknowledged only ten years after his death. These memoirs that he left, despite the fact that they all have the same aim, namely, to convince the reader that hyperbolic geometry exists, contain at some places remarks, examples, details and new computations which make them different from each other.
These memoirs include the \emph{New Elements of Geometry} \cite{Loba-New} (1835) (a revised version of the 1829  \emph{Elements of Geometry}), the \emph{Applications of Imaginary Geometry to Certain Integrals} \cite{Loba-Applications} (1836), \emph{Imaginary Geometry and its Applications} \cite{Loba-Imaginary-Russian}  (1835),  \emph{Imaginary Geometry} \cite{Loba-Imaginaire} 1837 (a revised version, in French,\footnote{Lobachevsky usually wrote in Russian, and he was a promoter of the use of Russian in scientific writing. This was uncommon in nineteenth-century Russia, where French was the official language for science, especially in mathematics. Russian was considered as inadequate. Lobachevsky nevertheless wrote his \emph{Imaginary Geometry}  and his \emph{Pangeometry} \cite{L} in French, and his  \emph{Geometrical Researches on the Theory of Parallels} \cite{Loba-Geometrische} in German, with the hope of being read in the West.} of the preceding memoir), the \emph{Geometrical Researches on the Theory of Parallels} \cite{Loba-Geometrische} (1840), in German, and finally, his \emph{Pangeometry} \cite{L} (1855), published first in Russian and then in French. The \emph{Pangeometry} is Lobachevsky's last work, written the year before his death. This memoir is a r\'esum\'e of Lobachevsky's work on non-Euclidean geometry and its applications, and it is probably his clearest account of the subject. It is also the conclusion of his lifetime work, and the last attempt he made to acquire recognition.\footnote{A recent edition of the \emph{Pangeometry}, translated into English and commented by the author of the present article, published by the European Mathematical Society Publishing House, is available, see \cite{L}. This edition also contains a complete and commented list of the various editions of Lobachevsky's writings.}

 The subject of the present paper is not the birth of hyperbolic geometry.\footnote{There are several good papers and books dealing with this subject, see e.g. \cite{Bonola}, \cite{Gray}, \cite{Gray2}, \cite{Greenberg}, \cite{Pont}, \cite{Rosenfeld} and the commentary in \cite{L}.} We shall rather focus on some important ideas contained in the \emph{Pangeometry} (and in the previous papers), which include the following:
\begin{enumerate}
\item The trigonometric formulae, which express the dependence of angles and edge lengths of triangles, are not only tools, but they are the basic elements of any geometry. In fact, Lobachevsky shows that all the analytical and geometrical theorems follow from these formulae.
\item Analysis in the sense of differential and integral calculus can be done in hyperbolic space without the use of any Euclidean model of hyperbolic space. (And indeed, Lobachevsky developed differential and integral calculus without knowing of any Euclidean model of his geometry. We recall by the way that the first such model was given by Beltrami in 1868, that is, more than twelve years after Lobachevsky's death).
\item There exist models of spherical and of Euclidean geometry within hyperbolic geometry.
\item If hyperbolic geometry were contradictory, then one could find a contradiction within Euclidean or spherical geometry.
\end{enumerate}
The fact that Lobachevsky expressed these ideas is rather unknown to mathematicians.
In the rest of this paper, we shall see how these ideas arise in Lobachevsky's work, and we shall comment on them.

\section{The spherical trigonometric formulae and the sphere model in hyperbolic space}\label{s:spherical}

Trigonometric formulae for spherical geometry, that is, equations giving relations between the side lengths and the angles at the vertices of spherical triangles, were known since Greek antiquity.\footnote{We recall that the Greeks did not use the sine function but a function they called \emph{chord}. The two functions are closely related since the sine of an arc is half of the double of its chord.} The classical proofs of these formulae use the embedding of the sphere in the ambient three-space, where the sides of triangles are arcs of great circles, that is, intersections of the sphere with planes passing through the origin, and the angles of triangles are the dihedral angles between the corresponding Euclidean planes.\footnote{The most important Greek spherical geometry treatise, which contains spherical trigonometry formulae,  is certainly Menelaus' \emph{Spherics} (second century A. D.) which reached us only through (several) Arabic commented translations, worked out between the eleventh and the thirteenth centuries. See e.g. \cite{RP1}, \cite{RP2} and \cite{RP3}.}  These formulae were recovered and reproved, using various methods, by many mathematicians over several centuries.\footnote{We mention by the way that Euler introduced an intrinsic method of proof of the spherical trigonometric formulae that uses the calculus of variations, see \cite{Euler}.} 

Lobachevsky, in his memoirs on geometry, established a complete set of spherical trigonometric formulae. His methods of proof are based on a construction in hyperbolic 3-space which uses the notion of parallelism in hyperbolic geometry. We recall this notion before reporting on that proof.

Let us first remind the reader that the notion of parallelism in non-Euclidean geometry is not identical to the Euclidean one.  in Euclidean geometry, two lines are said to be parallel if they are coplanar and if they do not intersect. Furthermore, in this setting, given a line $L$ and a point $x$ not on $L$, there is a unique line containing $x$ and disjoint from $L$, namely, the line containing $x$ and making a right angle with the segment joining perpendicularly $x$ to $L$.
     \begin{figure}[!hbp]
 \psfrag{x}{$x$}
\psfrag{L}{$L$}
\psfrag{I}{$I$}
\psfrag{N}{$N$}
\centering
\includegraphics[width=.6\linewidth]{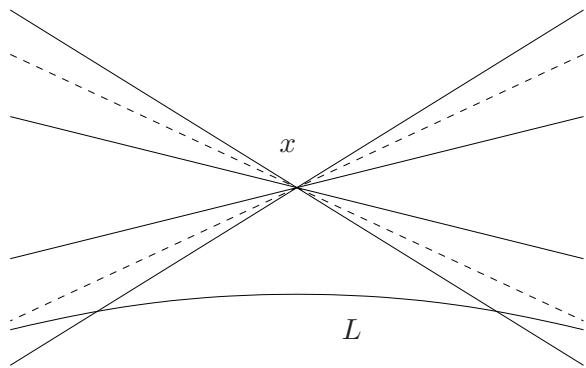}
\caption{\small{The two shaded lines are the parallels through $x$ to the line $L$.}}
\label{loba38}
\end{figure}

Now we consider the situation in the hyperbolic plane.  Let $L$ be a line and $x$ a point not on $L$. The line that contains $x$ and makes a right angle with the perpendicular from $x$ to $L$ does not intersect $L$ (otherwise we would have a triangle whose angle sum is greater than two right angles, which is not possible in neutral geometry, that is, in the geometry where one takes the Euclidean axioms without the parallel postulate, that is, this axiom is neutralized: it may hold or may not hold). However, this line is not the unique line that contains $x$ and that is disjoint from $L$. In fact, in hyperbolic geometry, one makes a distinction between two classes of disjoint pairs of lines, namely, {\it parallel} lines and {\it hyperparallel} lines.  A line {\it parallel} to $L$  through $x$ is a line that does not intersect $L$ and that is the limit of two families of lines, those that intersect $L$ and those that do not intersect $L$.\footnote{The definition of parallel line as a line separating the family of intersecting lines from the family of non-intersecting lines was used by the three founders of hyperbolic geometry.} With this definition, there are exactly two distinct lines through $x$ that are parallel to $L$ (see Figure \ref{loba38}). A line is said to be {\it hyperparallel} to $L$ if it is not parallel to $L$ and if it is disjoint from $L$. 

 In the memoirs that Lobachevsky wrote before the {\it Pangeometry}, e.g. in his  {\it Geometrische Untersuchungen zur Theorie der Parallellinien} to which he constantly refers in his {\it Pangeometry}, he studied in detail this relation of parallelism, showing in particular that it is symmetric and transitive.\footnote{In using the transitivity property, it is necessary to keep track of the direction of parallelism for pairs of parallel lines involved.}
 
 The definition of parallelism in dimension $3$ is based on the planar definition. More precisely, two lines in hyperbolic $3$-space are said to be parallel if they are coplanar and if they are parallel in the hyperbolic plane that contains them.

Now we can review Lobachevsky's proof of the spherical trigonometric formulae.

\begin{figure}[!hbp]
\psfrag{A}{$A$}
\psfrag{a}{$a$}
\psfrag{B}{$B$}
\psfrag{b}{$b$}
\psfrag{C}{$C$}
\psfrag{c}{$c$}
\psfrag{M}{$A'$}
\psfrag{N}{$B'$}
\psfrag{O}{$C'$}
\psfrag{m}{$m$}
\psfrag{n}{$n$}
\psfrag{k}{$k$}
\centering
\includegraphics[width=.6\linewidth]{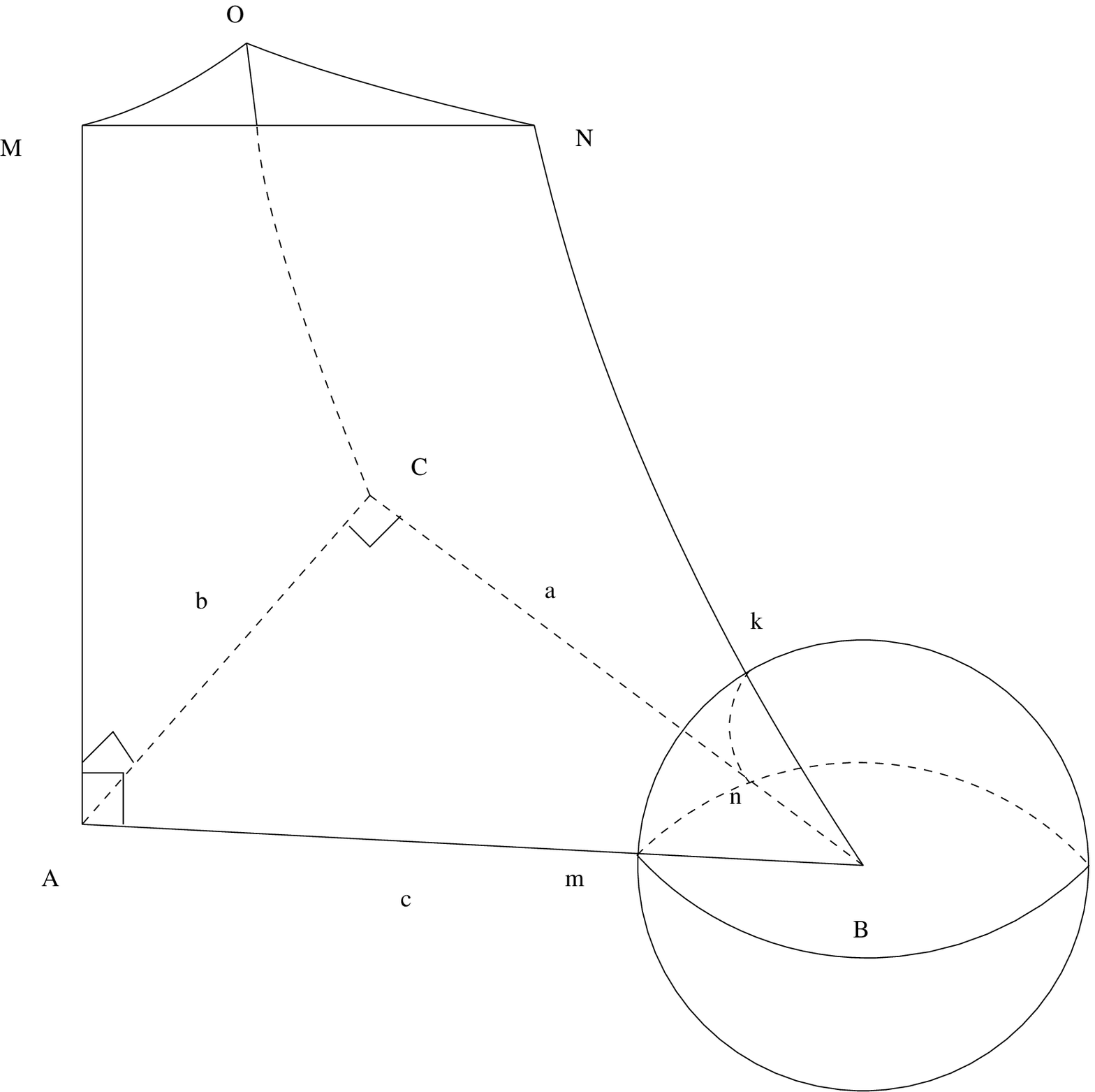}
\caption{\small{}}
\label{loba35}
\end{figure}

Lobachevsky starts with a hyperbolic triangle $ABC$ with a right angle at $C$ and with edges $a,b,c$ opposite respectively to the vertices $A,B,C$. He considers the plane containing this triangle as sitting in hyperbolic 3-dimensional space and he makes the following construction. Let $AA'$ be a geodesic line perpendicular to the plane of the triangle and consider the two planes  $BAA'$ and $CAA'$ (Figure \ref{loba35}). In the plane $BAA'$, construct the line $BB'$ parallel to the line $AA'$. Consider a third plane, containing the line $BB'$ and the point $C$. This plane intersects the plane $CAA'$ in a line $CC'$ parallel to $AA'$. (Lobachevsky proves this fact.)

Lobachevsky shows that the dihedral angles of the solid figure that is constructed over the base $ABC$ are completely determined by the angles of the triangle $ABC$.

Now consider in the 3-dimensional space a geometric sphere\footnote{A \emph{geometric sphere} in a metric space is the set of points that are at a certain distance from a given point, called the center of the sphere. Let us note for the attention of the non-geometrically knowledgeable reader that in a general metric space a geometric sphere may have properties which are very different from those of the usual sphere in Euclidean space; of course, it could be non-homeomorphic to such a sphere.} centered at $B$ and of radius smaller than $a$. It intersects the three planes passing by $B$ at three points $k,m,n$. 

In this situation, and in analogy with usual spherical geometry (the geometry of a sphere embedded in Euclidean 3-space), a spherical segment is defined as an arc in the intersection of the sphere by a (hyperbolic) plane passing through the center of the sphere, and a spherical triangle is defined as a triple of points joined by three spherical segments. (Strictly speaking, one has to consider only spherical triangles contained in a  hemisphere, in order to avoid pairs of points which can be joined by two distinct segments.) Distances and angles on this geometric sphere are also defined in analogy with those of spherical triangles on a standard sphere embedded in  Euclidean 3-space: the distance between two points on the sphere is the (hyperbolic) angular distance of the rays starting from the origin and passing through these points, and the angle between two lines is defined as the dihedral angle made by the planes that pass through the center of the sphere and that contain these two lines. For instance, in Figure \ref{loba35}, the angle at $k$ of the spherical triangle $kmn$ is the dihedral angle made by the planes $CBB'$ and $ABB'$, etc.

From the dihedral angles of the solid figure constructed on the base $ABC$, the angles of the spherical triangle $kmn$ are determined. Lobachevsky shows in particular that the angle at $m$ is right. He then deduces the trigonometric formulae for the right spherical triangle $kmn$ using the ambient hyperbolic space, as in the classical way where one deduces the trigonometric formulae for a spherical triangle using the embedding of the sphere in a Euclidean 3-space. 
 
 Passing to more standard notation where we have a spherical triangle $ABC$ with right angle at $C$ and with edges $a,b,c$ opposite to the vertices $A,B,C$ respectively, Lobachevsky writes the following formulae:

 \begin{equation}\label{eqn:spherical1}
   \left\{   \begin{array}{llll}
\sin A \sin c=\sin a\\
\cos b \sin A =  \cos B\\
\cos a \cos b= \cos c.
      \end{array}
          \right.
\end{equation}

\begin{figure}[!hbp]
 \psfrag{A}{$A$}
\psfrag{B}{$B$}
\psfrag{C}{$C$}
\psfrag{a}{$a$}
\psfrag{b}{$b$}
\psfrag{c}{$c$}
\centering
\includegraphics[width=.25\linewidth]{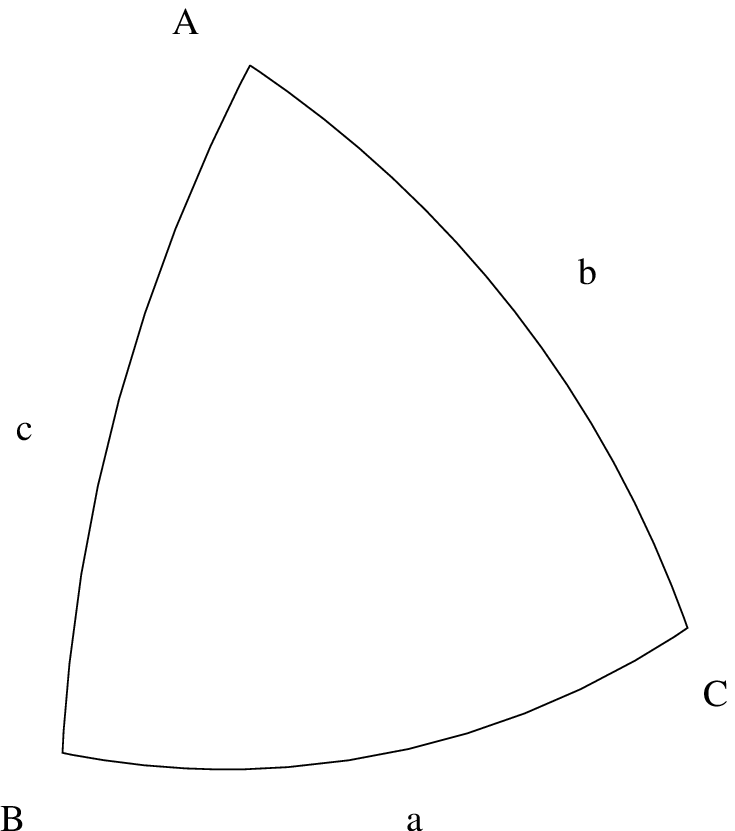}
\caption{\small{}}
\label{loba36}
\end{figure}

After these trigonometric formulae for the right triangle, Lobachevsky obtains the following formulae that are valid for an arbitrary spherical triangle, with sides $a,b,c$ opposite to angles $A,B,C$  (Figure \ref{loba36}):

 \begin{equation}\label{eqn:spherical1}
   \left\{   \begin{array}{llll}
\sin a \sin B=\sin b\sin A\\
\cos b-\cos a \cos c = \sin a \sin c \cos B\\
\cot a \sin b=\cot A\sin C+\cos b\cos C\\
\cos a\sin B\sin C=\cos B\cos C+\cos A.
      \end{array}
          \right.
\end{equation}

He then notices that these formulae for the sphere in hyperbolic space coincide with the formulae of usual spherical geometry. He concludes with the following (\cite{L} p. 22 of the English translation): 

\begin{quote}\small

It follows that spherical trigonometry stays the same, whether we adopt the hypothesis that the sum of the three angles of any rectilinear triangle is equal to two right angles, or whether we adopt the converse hypothesis, that is, that this sum is always less than two right angles.\footnote{At the end of \S 35 of his  {\it Geometrische Untersuchungen}, Lobachevsky drew
 a similar conclusion  (\cite{Loba-Geometrische}, translation in \cite{Bonola}): ``Hence spherical trigonometry is not dependent upon whether in a rectilinear triangle the sum of the three angles is equal to two right angles or not").} 
\end{quote}

We reformulate this statement as a theorem:

\begin{theorem}[Lobachevsky] \label{th:sphere}
The geometric spheres in hyperbolic space are models of spherical geometry.
\end{theorem}

This theorem is important and it was highly non-trivial a the time of Lobachevsky. Today, one can recover this theorem by arguing that a geometric sphere in three-dimensional hyperbolic space is a surface of constant intrinsic curvature, but it was not possible to make such a reasoning at Lobachevsky's time; this argument came later on, with the work of Riemann. As a matter of fact, Beltrami, as a conclusion of his 1868 paper \cite{Beltrami-Teoria}, obtained the same result, and he stated it in different words, using the concept of curvature. He writes (p. 62 of Stillwell's translation \cite{Sources}):
\begin{quote}
One sees that the geometry of spaces of constant positive curvature (which can appropriately be called \emph{spherical geometry} in the broad sense, since as equation (22) shows, the geodesic triangles are subject to the laws of spherical trigonometry) differs very markedly from \emph{pseudospherical} geometry,\footnote{Pseudo-spherical geometry is a word that Beltrami used for hyperbolic geometry.} even though both admit congruent figures. Moreover, pseudospherical geometry leads spontaneously to the consideration of spaces of positive curvature. In fact, if we put
\[\frac{a}{x}= y, \frac{x_1}{x}= y_1, ...,\frac{x_n}{x}= y_n,\]
in (26), we find
\[ds = R \sqrt{dy^2+dy_1^2+...+dy_n^2},\]
with the condition
\[y^2+y_1^2+...+y_n^2=1,\]
which, when we compare with equation (18) and take $\rho=$ const., tells us that the geodesic spheres of radius  $\rho$, in an $n$-dimensional space of constant curvature $-\frac{1}{R^2}$, are the $(n-1)$-dimensional spaces of constant curvature $\left( \frac{1}{R\sinh \frac{\rho}{R}}\right)^2.$ Thus, spherical geometry can be regarded as part of pseudospherical geometry.\end{quote}
\section{The horosphere as a model of Euclidean space}

Let us now state another result due to Lobachevsky.
\begin{theorem}[Lobachevsky] \label{th:eucli}
The horospheres in hyperbolic space are models of Euclidean geometry.
\end{theorem}

This is an interesting instance of a plane embedded in hyperbolic $3$-space which is naturally equipped with a Euclidean metric.

This result is stated in the first pages of the \emph{Pangeometry} (and also in Lobachevsky's previous memoirs). We present it in this section, before we proceed with the hyperbolic trigonometric formulae. Let us quote Lobachevsky. In the \emph{Pangeometry}, horocycles (respectively horospheres) are called limit circles (respectively horospheres).\footnote {In his {\it Geometrische Untersuchungen}  \S 34, Lobachevsky uses the word ``orisphere", which was transformed later on into ``horosphere".} 
After introducing the notion of horocycle with a given axis and passing through a certain point and after giving a method for constructing it, Lobachevsky writes the following (\cite{L}, p. 8 of the English translation):

\begin{quote}\small
The revolution of the limit circle around one of its axes\footnote{Again, the hyperbolic plane is considered as being embedded in hyperbolic $3$-space.} produces a surface which I call a {\it limit sphere}\index{limit sphere}\index{horosphere}, a surface which consequently is the limit towards which the sphere approaches as its radius increases to infinity. We shall call the   axis of revolution, and, consequently, any line that is parallel to this axis of revolution, an {\it axis}\index{limit sphere!axis}\index{axis of a limit sphere} of the limit sphere. We shall call a {\it diametral plane}\index{diametral plane}  any plane that contains one or more of the axes of the limit sphere. The intersections of the limit sphere with its diametral planes are limit circles. A portion of the surface of the limit sphere that is bounded by three limit circle arcs is called a {\it limit sphere triangle}\index{limit sphere triangle}. The limit circle arcs are called the {\it edges}, and the dihedral angles between the planes of these arcs are called the {\it angles} of the limit sphere triangle.

Two straight lines that are parallel to a third one are parallel ({\it Geometrische Untersuchungen} \S 25). It follows that all the axes of the limit circle and of the limit sphere are mutually parallel.

 If three planes intersect pairwise in three parallel lines and if we consider the region  that is situated between these parallels, the sum of the three dihedral angles that these planes form is equal to two right angles ({\it Geometrische Untersuchungen} \S 28).
 
 It follows from this theorem that the angle sum of any limit sphere triangle is equal to two right angles and, consequently, everything that we prove in ordinary geometry concerning the proportionality of edges of rectilinear triangles can be proved in the same manner in the Pangeometry of limit sphere triangles, by only replacing the lines that are parallel to one of the edges of the rectilinear triangle by limit circle arcs drawn from the points on one of the edges of the limit sphere triangle, and all of them making the same angle with that edge.\footnote{Lobachevsky explains here how one can construct similar (i.e. homothetic) triangles on the limit sphere, in a way analogous to the construction of similar triangles in the Euclidean plane. We recall that in neutral geometry, making the assumption that there exists a pair of similar non-congruent triangles is equivalent to adding Euclid's parallel axiom.}
\end{quote}

In other words, Lobachevsky deduces that the geometry of the limit sphere is Euclidean from the property that in that  geometry, the angle sum of triangles is equal to two right angles. It is known indeed that this property is equivalent to Euclid's parallel postulate, and one also has to check that the other axioms of Euclidean geometry are satisfied; this can be done indeed.

As he did for Theorem \ref{th:sphere}, Beltrami also rediscovered Theorem \ref{th:eucli}, and he formulated it in terms of differential geometry, in his 1868 paper {\it Teoria fondamentale degli spazii di curvatura costante} \cite{Beltrami-Teoria}. In fact, he obtained a concrete description of limit spheres on his  pseudo-sphere model of the hyperbolic plane. Using today's language, Beltrami's result says that the {\it length structure} induced by the Riemannian metric of $3$-dimensional hyperbolic space on a  limit sphere is Euclidean, whereas Lobachevsky's result is stated at the axiomatic level, namely, Lobachevsky showed that the geometry of the limit sphere satisfies the Euclidean geometry axioms.
 In a letter  to Ho\"uel, dated 1st of April 1868, Beltrami made the relation between the two points of view. He wrote (\cite{Beltrami-Boi} pp. 87--88): 
 \begin{quote}\small 
 From the formula $ds=\mathrm{const} \cdot \sqrt{d\eta_1^2+d\eta_2^2}$ that I established on page 21 of my last Memoir,\footnote{Beltrami refers here to his  \emph{Teoria fondamentale degli spazii di curvatura costante} \cite{Beltrami-Teoria}.} we can  deduce (or rather, we can check, because this is already contained in Lobachevsky) that the geometry of the limit sphere is nothing else than that of the Euclidean plane. This is exactly what I meant when I said that the curvature of this surface is zero. In other words, I meant that all the metric properties of this surface are the same as those of the ordinary plane, because the linear elements in both surfaces are identical. [...] For me, I would say that the limit sphere is one of the forms under which the Euclidean plane exists in the non-Euclidean space, \emph{considering the Euclidean plane as being defined by the property of having zero curvature}.
 \footnote{The letter is in French: [De la formule $ds=\mathrm{const} \cdot \sqrt{d\eta_1^2+d\eta_2^2}$, que j'ai \'etablie \`a la page 21 de mon dernier M\'emoire on tire (ou plut\^ot on v\'erifie, car cela se trouve dans Lobatcheffsky) que la g\'eom\'etrie de la sph\`ere-limite n'est pas autre chose que celle du plan euclidien. En disant que la courbure de cette surface est nulle je n'ai pas voulu dire autre chose. En d'autres termes j'ai voulu dire que toutes les propri\'et\'es m\'etriques de cette surface sont les m\^emes que celles du plan ordinaire, \`a cause de l'identit\'e des \'el\'ements lin\'eaires chez l'une et chez l'autre. [...] Pour mon compte je dirais que la sph\`ere-limite est une des formes sous lesquelles le plan euclidien existe dans l'espace non-euclidien, \emph{en consid\'erant le plan euclidien comme d\'efini par la propri\'et\'e d'avoir sa courbure nulle}.]} 
 \end{quote}
 
 A similar description of limit spheres is contained in Beltrami's letter to Ho\"uel, dated 12 October 1869, \cite{Beltrami-Boi} p. 87--88.

\section{The hyperbolic trigonometric formulae}\label{s:tr}

Lobachevsky stated his hyperbolic trigonometry formulae in terms of the {\it angle of parallelism function}. This is 
\begin{figure}[!hbp]
\psfrag{p}{$p$}
\psfrag{P}{$\Pi(p)$}
\centering
\includegraphics[width=.6\linewidth]{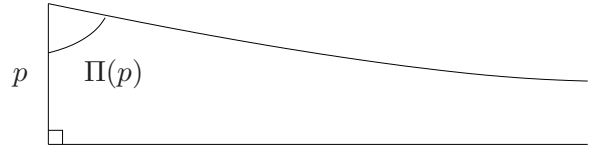}
\caption{\small{The parallelism angle $\Pi(p)$ corresponding to the segment $p$.}}
\label{parallelism}
\end{figure}
the function, denoted by $\Pi$,  which assigns to each segment of length $p$ the acute angle that this segment makes with a ray starting at an endpoint of this segment and which is parallel to a second ray that starts perpendicularly at the other endpoint of the segment (see Figure \ref{parallelism}). The value of the parallelism function depends only on the length $p$ of the segment and not on the segment itself.\footnote{This uses the fact that the hyperbolic plane is doubly homogeneous; that is, that any two segments of the same length are congruent (i.e. they are  equivalent through a motion of the plane). This fact is a consequence of the axioms of neutral geometry.} The parallelism angle of a segment of length $p$ is denoted by $\Pi(p)$.\footnote{In Lobachevsky's writings, as in other writings of the same period, and following the tradition of Euclid's {\it Elements}, the word ``line" (to which we have preferred the word ``segment" in the present translation) is often identified with the length of that line (i.e. of that segment). Thus, Lobachevsky says that $\Pi$ is a function of the line $p$, meaning that it is a function of the real number $p$.}

The \emph{angle of parallelism function}, or \emph{parallelism function},  $\Pi(p)$ was already introduced by Lobachevsky in his first written memoir, the  {\it Elements of geometry} (1829) \cite{Loba-Elements},  where it is denoted by $F(p)$. This function was used later on by several authors, with a reference to Lobachevsky. It is used in Beltrami's {\it Saggio di Interpretazione della geometria non-Euclidea} \cite{Beltrami-Saggio} and in Klein's {\it \"Uber die sogenannte Nicht-Euklidische Geometrie} \cite{Klein-Ueber}.

 The parallelism  function is related to the hyperbolic cosine function by the formula \[\sin \Pi(p)=\displaystyle\frac{1}{\cosh p},\] assuming the curvature of the space is $-1$.\footnote{The axioms of non-Euclidean geometry cannot specify the value of the curvature constant of the space. For a space of curvature $K=-{1}/{k^2}$, the formula becomes
  $\sin \Pi(p)= {1}/{\cosh (p/k)}$. We note by the way that the existence of this constant $k$ appears in the works of all three founders of hyperbolic geometry, although none of them did relate it to curvature. This constant appears in various forms, and we mention in this respect a letter from Gau\ss  \ to F. A. Taurinus, written on 8 November 1824: ``The assumption that in a triangle the sum of three angles is less than 180${}^{\mathrm{o}}$ leads to a curious geometry, quite different from ours, but thoroughly consistent, which I have developed to my entire satisfaction, so that I can solve every problem in it with the exception of the determination of a constant, which cannot be designated {\it a priori}. The greater one takes this constant, the nearer one comes to Euclidean geometry, and when it is chosen infinitely large, the two coincide."  (Greenberg's translation.)  [Die Annahme, dass die Summe der 3 Winkel kleiner sei als 180${}^{\mathrm{o}}$, f\"uhrt auf eine eigene, von der unsrigen (Euklidischen) ganz verschiedene Geometrie, die in sich selbt durchaus consequent ist, und die ich f\"ur mich selbst ganz befriedgend ausgebildet habe, so dass ich jede Aufgabe in derselben aufl\"osen kann mit Ausnahme der Bestimmung einer Constante, die sich apriori nicht ausmitteln l\"asst.]}

 \begin{equation}\label{eq:19}
\left\{
     \begin{array}{llll}
       \sin A \tan\Pi(a)=\sin B\tan\Pi(b)  \\
        1-\cos \Pi(b)\cos\Pi(c)\cos A=\displaystyle \frac{\sin\Pi(b)\sin\Pi(c)}{\sin \Pi(a)}   \\
          \cos A +\cos B\cos C=\displaystyle\frac{\sin B\sin C}{\sin \Pi(a)}\\
    \cot A \sin C\sin \Pi(b)+\cos C=\displaystyle\frac{\cos \Pi(b)}{\cos \Pi(a)}.\\
     \end{array}
  \right.
\end{equation}

\begin{figure}[!hbp]
\psfrag{A}{$A$}
\psfrag{a}{$a$}
\psfrag{B}{$B$}
\psfrag{b}{$b$}
\psfrag{C}{$C$}
\psfrag{c}{$c$}
\psfrag{M}{$A'$}
\psfrag{N}{$B'$}
\psfrag{O}{$C'$}
\psfrag{p}{$p$}
\psfrag{q}{$q$}
\psfrag{r}{$r$}
\centering
\includegraphics[width=.35\linewidth]{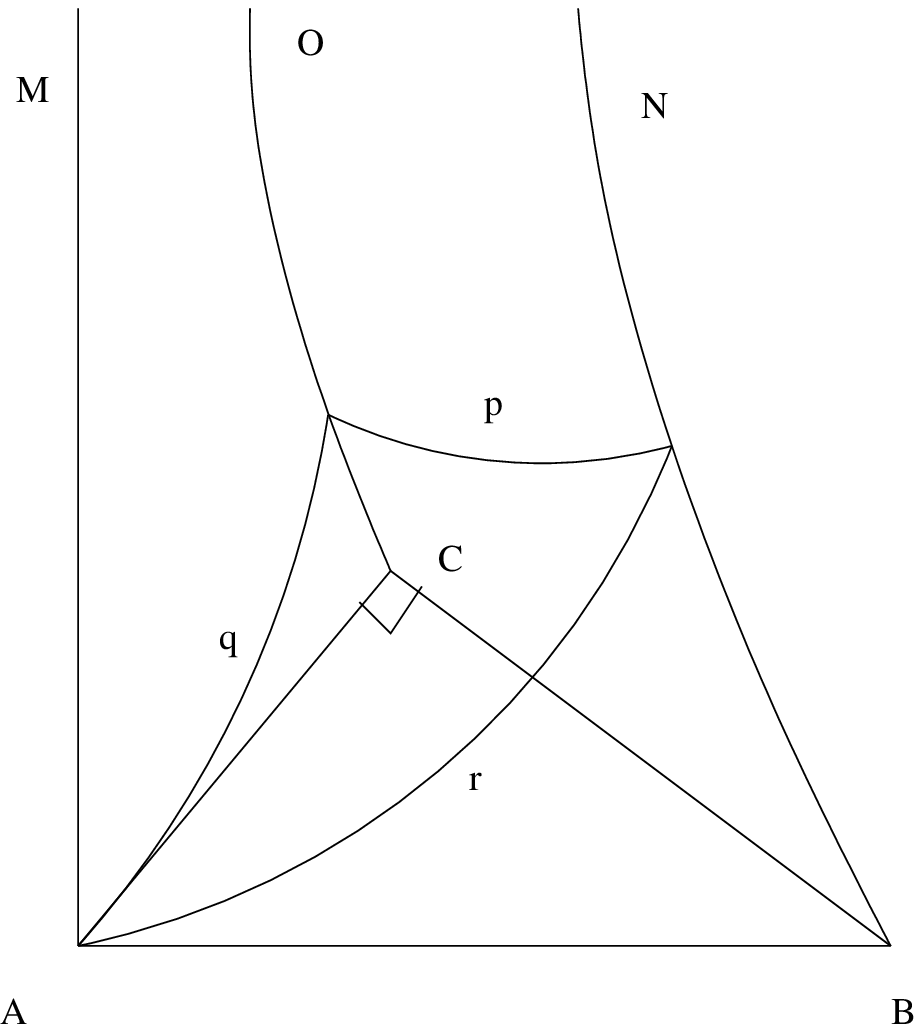}
\caption{\small{}}
\label{loba31}
\end{figure}

 To prove the trigonometric formulae for the right hyperbolic triangle $ABC$ (\S \ref{s:spherical} above), Lobachevsky constructs, besides the spherical triangle that we described there (Figure \ref{loba35}), the Euclidean triangle of Figure \ref{loba31}, that is, the triangle on the horosphere whose axis is the line $AA'$ and passing by $A$. This triangle is obtained as the intersection of that horosphere with the three planes that contain the edges of the triangle $ABC$ and the lines $AA', BB', CC'$.  Lobachevsky notes that this Euclidean triangle is again right, with right angle on $CC'$. Using the Euclidean trigonometric formulae for this triangle and the spherical trigonometric formulae he established for the triangle $kmn$, Lobachevsky obtains the trigonometric formulae for the hyperbolic triangle $ABC$. Then, decomposing again any triangle into two right triangle, he obtains a set of trigonometric formulae valid for any hyperbolic triangle $ABC$:

 We note that all these trigonometric formulae (spherical and hyperbolic) were already contained in Lobachevsky's first memoir, the {\it Elements of geometry}  (1828) \cite{Loba-Elements}, p. 21 of Engel's German translation \cite{Engel}, and also in his later works, such as the {\it Geometrische Untersuchungen} (\S 35 to 37) \cite{Loba-Geometrische}.

The rest of the \emph{Pangeometry}, where Lobachevsky develops the foundations of the analytic theory of hyperbolic geometry, is built on the hyperbolic trigonometric formulae. We shall not enter into the details here, but we can quote Lobachevsky:

 \begin{quote}\small Starting with these equations, Pangeometry becomes an analytic geometry, and thus it forms a complete and distinct geometric theory. Equations (\ref{eq:19}) are useful for representing curves by equations in terms of the coordinates of their points, and for calculating lengths and areas of curves, and areas and volumes of bodies, as I showed in the 1829 {\it Scientific Memoirs of the University of Kazan}.\footnote{Cf.  Lobachevsky's {\it Elements of geometry}, \cite{Loba-Elements}.}
 \end{quote}
 
\section{On the relative non-contradiction of hyperbolic geometry}
   
 The relative consistency (or non-contradiction) of non-Euclidean geometry with respect to the Euclidean one, at least in dimension two, is usually attributed to Beltrami, who established it in his famous paper {\it Saggio di Interpretazione della geometria non-Euclidea}  \cite{Beltrami-Saggio} by constructing a Euclidean model of it. In the introduction to that paper, explaining the goal of his paper, Beltrami writes the following:
 \begin{quote} \small We tried, to the limit of our capabilities,
 to make for ourselves an idea of the results to which the doctrine of Loba\-chevsky leads; and following a method which seems to us completely conformal to the good traditions of scientific investigation, we tried to provide that doctrine with a firm basis, before we admit for that theory the necessity of a new order of entities and of concepts.
 We believe that we succeeded in our goal with respect to the planar part of this doctrine, but we believe that this goal is impossible to attain for what concerns the rest.\footnote{The 3-dimensional case was dealt with soon after by Beltrami, in his paper \cite{Beltrami-Teoria}.}
 \end{quote}

It is interesting to note here a remark that Lobachevsky made, on the trigonometric formulae, at the end of his {\it Elements of geometry} (1829), which is related to the consistency issue of his newly discovered geometry. The remark is quoted in Rosenfeld \cite{Rosenfeld} p. 223 (Shenitzer's  translation from the Russian):
\begin{quote}
After we have found Equations (17) which represent the dependence of the angles and sides of a triangle; when, finally, we have given the general expressions for elements of lines, areas and volumes of solids, all else in the Geometry is a matter of analytics, where calculations must necessarily agree with each other, and we cannot discover anything new that is not included in these  first equations from which must be taken all relations of geometric magnitudes, one to another. Thus if one now needs to assume that some contradiction will force us subsequently to refute the principles that we accepted in this geometry, then such a contradiction can only hide in the very Equations (17). We note, however, that these equations become Equations (16) of spherical trigonometry as soon as, instead of the sides $a,b,c$ we put $a\sqrt{-1}, b\sqrt{-1},c\sqrt{-1}$; but in ordinary Geometry and in spherical Trigonometry there enter everywhere only ratios of lines; therefore ordinary Geometry, Trigonometry and the new Geometry will always agree among themselves.
\end{quote}

\begin{figure}[ht!]
\centering
 \psfrag{A}{\small $A$}
 \psfrag{B}{\small $B$}
 \psfrag{C}{\small $C$}
 \psfrag{a}{\small $a$}
 \psfrag{b}{\small $b$}
 \psfrag{c}{\small $c$}
 \psfrag{0}{\small $O$}
\includegraphics[width=.65\linewidth]{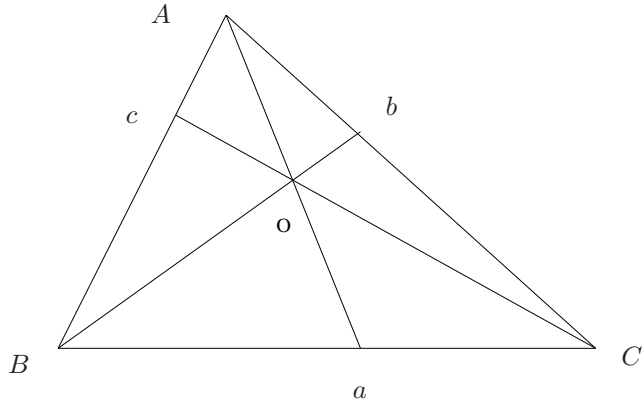}
\caption{\small {Equations (\ref{eq:Euler}) and (\ref{eq:Euler1})  give necessary and sufficient conditions for the lines  $Aa, Bb, Cc$ to have a common intersection point, in Euclidean and in spherical geometry respectively.}}
\label{concourantes}
\end{figure}

This is a first  indication of Lobachevsky's conviction that if there were a contradiction in hyperbolic geometry, then there would also be one in Euclidean geometry or in spherical geometry, and vice-versa. It is based on some formal analogy between the presentation of the formulae in the three geometries. Below, we shall point out another argument for the non-contradiction problem, which is also based on Lobachevsky's work. We note by the way that this formal analogy was noticed by several authors, before and after Lobachevsky.  For instance, J. H. Lambert (1728-1777), in his \emph{Theorie der Parallellinien}, written in 1766 and published posthumously, (see \cite{Engel-Staeckel} and the forthcoming edition \cite{Lambert}) observed that certain geometrical properties of hyperbolic geometry -- which for him was hypothetical\footnote{Lambert is among the immediate predecessors of hyperbolic geometry who developed a complete theory of that geometry with the hope of finding a contradiction.} -- are obtained from analogous properties of spherical geometry by multiplying some distances by the imaginary number $\sqrt{-1}$. He declared then that one should almost conclude that this geometry takes place on a ``sphere of imaginary radius".
 F. A. Taurinus  (1794--1874) (a contemporary of Lobachevsky) made a similar remark. In his memoir, also called {\it Theorie der Parallellinien} (1825) \cite{Engel-Staeckel}, Taurinus obtained fundamental trigonometric formulae for hyperbolic  geometry (which was, from his point of view, like it was for Lambert, purely hypothetical). He noticed the formal analogy between the spherical and the hyperbolic cases and he also declared that the hyperbolic trigonometry formulae are obtained by working on a sphere of ``imaginary radius". Likewise,  Beltrami, in his  {\it Saggio di interpretazione della geometria non-euclidea} \cite{Beltrami-Saggio}, showed that the trigonometric formulae for the pseudo-sphere (which is a model he had constructed for the hyperbolic plane) can be obtained from those of the usual sphere by considering the pseudo-sphere as a sphere of imaginary radius $\sqrt{-1}$, and he attributed this observation to E. F. A. Minding  and to D. Codazzi. Soon after Beltrami, Klein, in his \emph{On the so-called non-Euclidean geometry} (1871), while he worked for a common ground for the three geometries (Euclidean, spherical and hyperbolic) in the setting of projective geometry, made the same remark: ``The trigonometric formulae that hold for our measure result from the formulae of spherical trigonometry by replacing sides by sides divided by $\frac{c}{i}$." (\cite{Klein-Ueber} \S 12). We also refer the interested reader to the comments in the English translation of the \emph{Pangeometry} \cite{L}. 
 Let us also note that Euler had also noticed the analogies between Euclidean and spherical geometry formulae, and in several memoirs, he proved Euclidean theorems and then their spherical analogues. For instance, in his memoir \emph{Geometrica et sphaerica quaedam} (Concerning Geometry and Spheres)  \cite{Euler-Geometrica-T} Euler considers the following problem: Given a triangle and three segments, each segment joining the vertex of a triangle to the line containing the opposite side, find a condition so that the three segments lie on three lines that have a common point. Then, denoting by  $A,B,C$ the vertices of the triangle, $a,b,c$ the points on the respective sides and $O$ the intersection point of the three lines (Figure \ref{concourantes}), Euler finds, in the Euclidean case, the relation 
 \begin{equation}
\label{eq:Euler}
\frac{AO}{Oa}\times\frac{BO}{Ob}\times\frac{CO}{Oc}=\frac{AO}{Oa}+\frac{BO}{Ob}+\frac{CO}{Oc}+2.
\end{equation}
 In the spherical case, he finds the relation
  \begin{equation}
\label{eq:Euler1} \frac{\tan AO}{\tan Oa}\frac{\tan BO}{\tan Ob}\frac{\tan CO}{\tan Oc} =  \frac{\tan AO}{\tan Oa} + \frac{\tan BO}{\tan Ob} + \frac{\tan CO}{\tan Oc} + 2.
\end{equation}

In other words, to pass from the Euclidean to the spherical cases, one simply replaces side lengths by the tangents of these side lengths.

Thus if we believe, like Lobachevsky, that the trigonometric formulae are at the basis of all of geometry, then looking at the formal analogies between the formulae in the three geometries, one can try to argue (as Lobachevsky did) that if there were a contradiction in one of the geometries then there would also be one in the other two, or at least in one of the other two.
 
 There is another argument that Lobachevsky could have used for the consistency of his geometry. This is the fact that the \emph{proof} of the trigonometric formulae of hyperbolic geometry are derived from the formulae of Euclidean and of spherical geometry, as we pointed out in \S \ref{s:spherical} and \S \ref{s:tr} and, therefore, one might conclude that the (hypothetic) existence of a contradiction in hyperbolic geometry would come from a contradiction in Euclidean or in spherical geometry.

 Finally, let us recall that the consistency question was raised in precise terms by David Hilbert, who reduced the consistency of the axiom set of Euclidean geometry (and of hyperbolic and spherical
geometries) to that of arithmetic. We also recall that G\"odel,  in 1931, proved that arithmetic with multiplication cannot be shown to be consistent within arithmetic. In conclusion, the consistency question of non-Euclidean (and Euclidean !) geometries has been a long and difficult issue and it is not surprising that Lobachevsky, who constantly thought about the consistency of his geometry, did not succeed in settling it, although there is more than enough evidence in his works to conclude with the relative consistency.


\begin{thebibliography}{ABCD}
 

\bibitem{Beltrami-Saggio} E. Beltrami, Saggio di Interpretazione della geometria non-Euclidea. (Essay on the interpretation of non-Euclidean geometry) \textit{Giornale di Matematiche} vol. VI (1868).  Beltrami's {\it Works}, Vol. I, pp. 374--405. English translations by J. Stillwell,  in \cite{Sources} pp. 7--34. 


\bibitem{Beltrami-Teoria} E. Beltrami,  Teoria fondamentale degli spazii di curvatura costante. (Fundamental theory of spaces of constant curvature). \textit{Annali di matematica pura ed applicata} serie II, Tomo II (1868-69), pp. 232--255. Beltrami's {\it Works}, Vol. I, pp. 406--429. French translation, by J. Ho\"uel, \emph{Th\'eorie fondamentale des espaces de courbure constante},    \textit{Annales scientifiques de l'\'Ecole Normale Sup\'erieure} S\'er. 1, Tome VI,  (1869), pp.  347--375.   English translation by J. Stillwell,  in \cite{Sources} pp. 41--62. 



 \bibitem{Beltrami-Boi} L. Boi, L. Giacardi and R. Tazzioli (ed.) {\it La d\'ecouverte de la g\'eom\'etrie non euclidienne sur la pseudosph\`ere. Les lettres d'Eugenio Beltrami \`a Jules Ho\"uel (1868-1881).} (The discovery of non-Euclidean geometry on the pseudosphere. Eugenio Beltrami's letters to Jules Ho\"uel (1868-1881).)
 Introduction, notes and commentary by L. Boi, L. Giacardi and R. Tazzioli. Preface by Ch. Houzel and E. Knobloch. 
Collection Sciences dans l'Histoire. Paris, Librairie Scientifique et Technique Albert Blanchard (1998). 


\bibitem{Bonola} R. Bonola,
\textit{La geometria non-euclidea, Esposizione storico-critica del suo sviluppo}.  First edition, Ditta Nicola Zanchinelli editore, Bologna 1906.  English
  translation by H. S. Carslaw: \textit{Non-Euclidean Geometry, A critical and historical study of its development},1912.  Reprinted by Dover, 1955.


  \bibitem{Coxeter-Functions} H. S. M. Coxeter,    The functions of Schl\"afli and Lobatschefsky. \textit{Quart. J. Math.} 6 (1935), 13-29.  
  
   

 \bibitem{Engel} F. Engel, Nikolaj Iwanowitsch Lobatschefskij, Zwei Geometrische Abhandlungen aus dem russischen \"ubersetzt mit Ammerkungen und mit einer Biographie des Verfassers. German translation of Lobachevsky's \emph{Elements of Geometry} and \emph{New elements of geometry, with a complete theory of parallels}. Leipzig, Teubner, Leipzig 1898.
  
  
\bibitem{Euler} L. Euler, Principes de la trigonom\'etrie sph\'erique tir\'es de la m\'ethode des plus grands et des plus petits, \emph{M\'emoires de l'Acad\'emie des sciences de Berlin} 9 (1753), 1755, p. 233-257 ; \emph{Opera omnia}, Series 1, vol. XXVII, p. 277-308.



\bibitem{Euler-Pappi-T} L. Euler, Problematis cuiusdam Pappi Alexandrini constructio,  \emph{Acta academiae scientarum imperialis Petropolitinae} 4, 1783, p. 91-96 ; \emph{Opera omnia} Series 1, vol. XXVI, p. 237-242.

\bibitem{Euler-Geometrica-T}  L. Euler,  Geometrica et sphaerica quaedam, \emph{M\'emoires de l'Acad\'emie des sciences de Saint-P\'etersbourg} 5, 1815, p. 96-114 ; \emph{Opera omnia} Series 1, vol. XXVI, p. 344-358




\bibitem{Gauss} C. F. Gauss, \textit{Collected Works}, Vol. VIII,
K\"onigliche Gesellshaft der Wissenschaften, G\"ottingen, 1900.



\bibitem{Gray} J. J. Gray, \textit{Worlds Out of Nothing.
A Course in the History of Geometry in the 19th Century}.
 Springer Undergraduate Mathematics Series,
2007. Second edition, 2011.



\bibitem{Gray2} J. J. Gray, \textit{History of Geometry}, in: the Princeton Companion to Mathematics, Princeton University Press, 2006, pp. 83--95.


\bibitem{Greenberg} M. J. Greenberg, \textit{Euclidean and Non-Euclidean Geometries}. 4th edition, New York, N.Y., W. H. Freeman, 2008  (first edition 1980).


\bibitem{Klein-Ueber} F. Klein,  \"Uber die sogenannte Nicht-Euklidische Geometrie (erster Aufsatz). (On the so-called non-Euclidean geometry, Part I) \textit{ Mathematische Annalen} IV (1871), 573--625.  English translation by J. Stillwell in \cite{Sources}.

\bibitem{Lambert} J. H. Lambert,  \emph{Theorie der Parallellinien}, in \cite{Engel-Staeckel}. Critical edition with a French translation by A. Papadopoulos and  G. Th\'eret, to appear, \'Ed. Albert Blanchard, Paris, Collection Sciences dans l'Histoire,  2014.

 \bibitem{Loba-Elements} N. I. Lobachevsky, Elements of geometry (Russian),  \textit{Kazansky Vestnik},  Issue 25: Feb. and Mar., 1829, pp. 178--187; Apr., 1829 pp. 228--241. Issue 27: Nov. and Dec., 1829, pp. 227--243. Issue 28: Mar. and Apr., 1830, pp. 251--283; Jul. and Aug., 1830, pp. 571--683.    German translation in \cite{Engel}. 


\bibitem{Loba-Applications} N. I. Lobachevsky,   Application of Imaginary Geometry to Certain Integrals (Russian)
\textit{Uchenye zapiski kazanskogo imperatorskogo  universiteta}. 1836. Issue I. pp. 3--166.






 \bibitem{Loba-Geometrische} N. I. Lobachevsky,
 \emph{Geometrische Untersuchungen zur theorie der Parallellinien} (Geometrical Researches on the Theory of Parallels). Kaiserl. russ. wirkl. Staatsrathe und ord. Prof. der Mathematik bei der Universit\"at Kasan.  Berlin 1840. In der Finckeschen Buchhandlung. 61 pages. 
 English translation by G. B. Halsted, first published by the University of Texas at Austin, 1891. Reprinted in Bonola \cite{Bonola}.   


 \bibitem{Loba-Imaginary-Russian}    N. I. Lobachevsky,   Imaginary geometry (Russian)  \textit{Uchenye zapiski kazanskogo imperatorskogo  universiteta}. 1835, Issue I, pp. 3-88. 
 
 \bibitem{Loba-New}   N. I. Lobachevsky, \emph{New elements of geometry, with a complete theory of parallels} (Russian) Uchenye zapiski kazanskogo imperatorskogo  universiteta 1835, Issue 3, pp. 3--48; 1836, Issue 2, pp. 3--98 \& Issue 3, pp. 3--50; 1837, issue 1, pp. 3--97; 1838, Issue 1, pp. 3--124 \& Issue 3, pp. 3--65.
   German translation of Chapters I to XI by F. Engel in \cite{Engel}, pp.  67--236.
 French translation, by  F. Mailleux,  Nouveaux principes de le g\'eom\'etrie avec une th\'eorie compl\`ete des parall\`eles. \textit{M\'em. Soc. Royale sci. de Li\`ege} (3), 2 No. 5, pp. 101 ; 3--101 No. 2, pp. 32--32 (1899). 



\bibitem{Loba-Imaginaire}  N. I. Lobachevsky,   G\'eom\'etrie imaginaire, \textit{J. Reine Angew. Math}. 17, 295-320 (1837). 




   \bibitem{L} N. I.  Lobachevsky, Pangeometry, Russian version in  \textit{Uchenye zapiski kazanskogo imperatorskogo  universiteta}. 1855, Issue 1, pp. 1--76.  French version, \emph{Pang\'eom\'etrie ou pr\'ecis de g\'eom\'etrie fond\'ee sur une th\'eorie g\'en\'erale et rigoureuse des parall\`eles}.  in: ``Recueil d'articles \'ecrits par les professeurs de l'Universit\'e de Kazan \`a l'occasion du cinquantenaire de sa fondation" (Collection of memoirs written by professors of the University of Kazan on the occasion of the 50th anniversary of its foundation), Vol. I, pp.  277--340 (1856). New edition with an English translation and a commentary by A. Papadopoulos, 
 Heritage of European Mathematics, Vol. 4, European Mathematical Society, Z\"urich, 2010.
 

\bibitem{Milnor} J. Milnor, Hyperbolic geometry: The first 150 years. \textit{Bull. Amer. Math. Soc.} 6 (1982), 9--24.


 \bibitem{Pont} J.-C. Pont,  {\it L'Aventure des parall\`eles  : Histoire de la g\'eom\'etrie non euclidienne: Pr\'ecurseurs et attard\'es}. Peter Lang ed. 1986.

 \bibitem{RP1} R. Rashed and A. Papadopoulos, On Menelaus' \emph{Spherics} III.5 in Arabic mathematics, I: Ibn `Ir\=aq,  \emph{Arabic Science and Philosophy},  vol. 24 (2014), pp. 1--68 (Cambridge University Press).

 
  
 
  \bibitem{RP2}  R. Rashed and A. Papadopoulos, On Menelaus' \emph{Spherics} III.5 in Arabic mathematics,  II: Al-Haraw\=\i , Na\d{s}\=\i r al-D\=\i n al-\d{T}\=us\=\i  \  and Ibn Ab\=\i\  Jarr\=ada, \emph{Arabic Sciences and Philosophy}, Cambridge University Press, to appear in 2014.
 

  
   \bibitem{RP3} R. Rashed and A. Papadopoulos, \emph{Menelaus' \emph{Spherics}}, Critical edition with an English translation, to appear.
   

\bibitem{Rosenfeld} B. A. Rosenfeld,  {\it History of non-Euclidean geometry}, translated by Abe Shenitzer, Studies in the History of Mathematics and Physical Sciences 12, Springer Verlag, 1988.




\bibitem{Sources} J. Stillwell (editor), {\it Sources of Hyperbolic geometry}, History of Mathematics Vol. 10, AMS-LMS,  1996.


 \bibitem{Engel-Staeckel}  P.~St\"ackel  et F.~Engel,  \emph{Die Theorie der Parallellinien von Euklid bis auf Gauss, eine Urkundensammlung zur Vorgeschichte der nicht-euklidischen Geometrie},  B. G. Teubner, Leipzig, 1895.

\bibitem{Thurston} W. P. Thurston, \emph{The
Geometry and Topology of Three-manifolds}, Mimeographed notes, Princeton University, 1979.

\end{thebibliography}
\end{document}